# Limit distributions for the problem of collecting pairs

PAVLE MLADENOVIĆ

*University of Belgrade, Faculty of Mathematics, Studentski trg 16, 11000 Belgrade, Serbia.
E-mail: paja@matf.bg.ac.yu*

Let $N_n = \{1, 2, \ldots, n\}$. Elements are drawn from the set $N_n$ with replacement, assuming that each element has probability $1/n$ of being drawn. We determine the limiting distributions for the waiting time until the given portion of pairs $jj$, $j \in N_n$, is sampled. Exact distributions of some related random variables and their characteristics are also obtained.

## 1. Introduction

Combinatorial problems in the theory of probability and mathematical statistics have been studied extensively. Many of them are formulated in the form of urn models. In such problems, one usually considers a sequence of experiments with some stopping rule defined a priori and the problem is to determine the exact and/or limit distribution of the waiting time until the last experiment. Sums and extreme values of random variables and rare events in a sequence of experiments appear naturally in connection with problems of this kind. Consequently, many different approaches, methods and techniques have been used to investigate combinatorial problems from the probabilistic point of view and a number of limit theorems have been proven. The method of characteristic and moment generating functions in summing random variables was used by Erdös and Rényi (1961), Békéssy (1964), Baum and Bilingsley (1965), Holst (1971), Samuel-Cahn (1974) and Flato (1982). The method of embedding in Poisson processes was used by Holst (1977, 1986). For a general list of references concerned this subject, see, for example, Johnson and Kotz (1977), Kolchin, Sevastyanov and Chistyakov (1976), Kolchin (1984) and Barbour, Holst and Janson (1992).

In this paper, the following problem will be studied. We sample with replacement from the set $N_n = \{1, 2, \ldots, n\}$, under the assumption that each element of $N_n$ has probability $1/n$ of being drawn, and we are interested in the waiting time until a given portion of pairs $jj$, $j \in N_n$, is sampled. In order to get limit distributions, we shall use the method of characteristic functions and also an approach based on the extreme value theory for stationary sequences; see Leadbetter, Lindgren and Rootzén (1983). The problem we are going to consider is a generalization of the coupon collector's problem. Originally, the waiting time for all $j$'s from $N_n$, supposing that all elements from $N_n$ have equal probability to be drawn at each step, was named the *coupon collector's problem*. The limiting distribution for this problem was first determined by Erdös and Rényi (1961).

A natural generalization of the coupon collector's problem is the problem of possible limiting distributions for the waiting time for a given portion of $j$'s from $N_n$. This problem was solved by Baum and Bilingsley (1965). Another generalization is the waiting time problem for a given number of appearances of all $j$'s from $N_n$. The limiting distribution for this problem and some related results were obtained by different authors employing different methods; see, for example, Holst (1986) and Mladenović (1999, 2006).

The problem of waiting time until a given portion of pairs $jj$, $j \in N_n$, is sampled can also be considered by using an approach based on point process theory. For a presentation of this theory, see, for example, Chapter 3 of Resnick (1987). Let $A_n(k)$ be the number of different pairs $jj$, $j \in N_n$, sampled until to the $k$th drawing and let $\{L_n(k), k \geq 1\}$ be the point process determined by the indices where the process $A_n$ jumps. Point process theory enables analysis of the asymptotic behavior of $\{L_n(k), k \geq 1\}$ and $\{A_n(L_n(k)), k \geq 1\}$. This approach was used in Chapter 4 of Resnick (1987) to study the structure and asymptotic behavior of records in a sequence of i.i.d. random variables with a continuous distribution function $F$. However, the underlying distributions in the problem that will be considered in this paper are discrete and depend on $n$.

The paper is organized as follows. Section 2 contains preliminaries, necessary notation and auxiliary results concerning exact distributions of random variables that appear in connection with the problem considered. Main results on asymptotic distributions are formulated in Section 3. Proofs of theorems from Sections 2 and 3 are given in Sections 4 and 5.

## 2. Preliminaries, notation and auxiliary results

Let $Z_1$, $Z_2$, $Z_3$, ... be a sequence of independent random variables with the uniform distribution over the set $N_n = \{1, 2, \ldots, n\}$. Throughout this paper, we shall use the following notation:

$$X_{nj} = \min\{k : Z_{k-1} = Z_k = j\}, \qquad j \in N_n \text{ is a fixed number}, \qquad (2.1)$$

$$\widetilde{Y}_{nj} = \min\{k : Z_{k-1} = Z_k = a \text{ for some } a \in A \subset N_n, |A| = j\}, \qquad (2.2)$$

$$M_n = \max\{X_{n1}, X_{n2}, \ldots, X_{nn}\}, \qquad (2.3)$$

$$M_n^{(k)} = \text{the } k\text{th maximum of random variables } X_{n1}, \ldots, X_{nn}, \qquad (2.4)$$

where $|A|$ is the number of elements of a set $A$. $X_{nj}$ is then the waiting time until the pair $jj$ for some fixed $j \in N_n$ occurs as a run in the process $Z_1, Z_2, \ldots, \widetilde{Y}_{nj}$ is the waiting time until some pair $aa$, where $a \in A$ and $|A| = j$, occurs as a run in the same process and $M_n$ is the waiting time until all $n$ pairs $jj$, $j \in N_n$, occur.

Let $Y_{nn}$ be the waiting time until the first pair $j_1 j_1$, where $j_1 \in N_n$, occurs as a run in the process $Z_1, Z_2, \ldots$. Let $Y_{n,n-1}$ be the waiting time for the second pair $j_2 j_2$, where $j_2 \in N_1 \setminus \{j_1\}$, after the occurrence of the first pair, etc. Then $Y_{nj} \stackrel{d}{=} \widetilde{Y}_{nj}$ for any $j \in \mathbf{N}_n$, where $X \stackrel{d}{=} Y$ means that random variables $X$ and $Y$ have the same distribution. Let us denote by $S_{n,a_n}$ the waiting time until $a_n$ of the pairs $jj$, $j \in N_n$, occur, that is,

$$S_{n,a_n} = Y_{nn} + Y_{n,n-1} + \cdots + Y_{n,n-a_n+1}, \qquad a_n \in N_n. \tag{2.5}$$

It is obvious that the following relations hold:

$$Y_{n1} \stackrel{d}{=} X_{n1} \stackrel{d}{=} X_{n2} \stackrel{d}{=} \cdots \stackrel{d}{=} X_{nn}; \tag{2.6}$$

$$S_{nn} = Y_{nn} + Y_{n,n-1} + \cdots + Y_{n1}$$
$$= \max\{X_{n1}, \ldots, X_{nn}\} = M_n; \tag{2.7}$$

$$S_{n,n-k+1} = M_n^{(k)}, \quad k \text{ is a fixed positive integer}; \tag{2.8}$$

$$P\{X_{n1} > m, X_{n2} > m, \ldots, X_{nj} > m\} = P\{Y_{nj} > m\}. \tag{2.9}$$

It is also clear that random variables $Y_{nn}, Y_{n,n-1}, \ldots, Y_{n1}$ are independent, but random variables $X_{n1}, X_{n2}, \ldots, X_{nn}$ are dependent. Let $F_n(x)$ be the common distribution function of random variables $X_{n1}, \ldots, X_{nn}$. As usual, $\Phi(x)$ is the standard normal distribution function. First, we shall give exact distributions and related characteristics of random variables $X_{nj}$ and $Y_{nj}$ and results concerning the asymptotic behavior of mean and variance of the random variable $S_{n,a_n}$.

**Theorem 2.1.** (a) *The distribution of the random variable $X_{nj}$ is given by*

$$P\{X_{nj} = k\} = \sum_{s=0}^{[k/2]-1} \binom{k-s-2}{s} \left(1 - \frac{1}{n}\right)^{k-s-2} \frac{1}{n^{s+2}}, \qquad k \geq 2. \tag{2.10}$$

(b) *If $u_n = n^2(x + \ln n)$, then the following equality holds:*

$$\lim_{n \to \infty} n(1 - F_n(u_n)) = \lim_{n \to \infty} nP\{X_{nj} > u_n\} = e^{-x}. \tag{2.11}$$

**Theorem 2.2.** (a) *The distribution of random variable $Y_{nj}$ is given by*

$$P\{Y_{nj} = k\} = \frac{j}{n\{(n+1)^2 - 4j\}^{1/2}} \cdot \left\{ \left(\frac{t_1}{n}\right)^{k-1} - \left(\frac{t_2}{n}\right)^{k-1} \right\}, \qquad k \geq 2, \tag{2.12}$$

where

$$t_1 = t_1(j) = \frac{n - 1 + \{(n+1)^2 - 4j\}^{1/2}}{2}, \tag{2.13}$$

$$t_2 = t_2(j) = \frac{n - 1 - \{(n+1)^2 - 4j\}^{1/2}}{2}. \tag{2.14}$$

(b) *Exact values of the mean and variance of the random variable* $Y_{nj}$ *are given by*

$$EY_{nj} = \frac{n^2 + n}{j}, \qquad \operatorname{var} Y_{nj} = \frac{n^4}{j^2}\left(1 + \frac{2}{n} - \frac{3j-1}{n^2} - \frac{j}{n^3}\right). \tag{2.15}$$

**Theorem 2.3.** *The asymptotic behavior of the mean* $\mu_n$ *and the variance* $\sigma_n^2$ *of the random variable* $S_{n,a_n}$ *is determined as follows:*

(a) *if* $a_n \to \infty$ *and* $a_n/n \to 0$ *as* $n \to \infty$, *then*

$$\mu_n = -n^2 \ln\left(1 - \frac{a_n}{n}\right) + o(na_n^{1/2}), \qquad \sigma_n^2 \sim n^2 a_n \qquad \text{as } n \to \infty; \tag{2.16}$$

(b) *if* $a_n/n \to \lambda \in (0,1)$ *as* $n \to \infty$, *and* $\lambda_0 = \lambda/(1-\lambda)$, *then*

$$\mu_n = -n^2 \ln\left(1 - \frac{a_n}{n}\right) + o(n^{3/2}), \qquad \sigma_n^2 \sim \lambda_0 n^3 \qquad \text{as } n \to \infty; \tag{2.17}$$

(c) *if* $a_n/n \to 1$ *and* $b_n = n - a_n \to \infty$ *as* $n \to \infty$, *then*

$$\mu_n = -n^2 \ln\left(1 - \frac{a_n}{n}\right) + o(n^2 b_n^{-1/2}), \qquad \sigma_n^2 \sim n^4/b_n \qquad \text{as } n \to \infty; \tag{2.18}$$

(d) *as* $n \to \infty$, *the mean and variance of the maximum* $M_n = S_{nn}$ *are given by*

$$EM_n = (n^2 + n)(\ln n + \gamma) + \frac{n}{2} + \frac{5}{12} - \frac{1}{12n} + \frac{1}{120n^2} + o\left(\frac{1}{n^2}\right), \tag{2.19}$$

$$\operatorname{var} M_n = \frac{\pi^2 n^4}{6} + \left(\frac{\pi^2}{3} - 1\right) n^3 - 3n^2 \ln n + O(n^2), \tag{2.20}$$

*where* $\gamma = 0.5772156649\ldots$ *is the Euler constant.*

## 3. Main results

The next three theorems give the limiting distribution of the random variable $S_{n,a_n}$ for different types of asymptotic behavior of the sequence $(a_n)$.

**Theorem 3.1.** *If* $a_n = k$ *for every* $n \in N$, *where* $k$ *is a fixed positive integer, then the random variable* $n^{-1} S_{n,a_n}$ *converges in distribution to a random variable whose characteristic function is*

$$f(t) = (1 + t^2)^{-k/2} e^{ik \cdot \arctan t}. \tag{3.1}$$

If $a_n/n \to \lambda \in [0,1]$, $a_n \to \infty$ and $b_n = n - a_n \to \infty$ as $n \to \infty$, then $(S_{n,a_n} - \mu_n)/\sigma_n$ has asymptotically normal $(0,1)$ distribution, where $\sigma_n^2 = \operatorname{var} S_{n,a_n}$ and $\mu_n = \mathrm{E} S_{n,a_n}$. Denote by $\nu_n$ and $\tau_n^2$ the main terms of $\mu_n$ and $\sigma_n^2$, respectively, that are determined by (2.16)–(2.18). Using the Khinchine lemma, we can conclude that the constants $\mu_n$ and $\sigma_n^2$ can be replaced in limit theorems by $\nu_n$ and $\tau_n^2$ because $\sigma_n/\tau_n \to 1$ and $(\mu_n - \nu_n)/\sigma_n \to 0$ as $n \to \infty$. More detailed results are provided by the following theorem.

**Theorem 3.2.** (a) *If $a_n \to \infty$ and $a_n/n \to 0$ as $n \to \infty$, then*

$$\lim_{n \to \infty} P\left\{ \frac{S_{n,a_n} + n^2 \ln(1 - a_n/n)}{na_n^{1/2}} \leq x \right\} = \Phi(x). \tag{3.2}$$

(b) *If $a_n/n \to \lambda \in (0,1)$ as $n \to \infty$ and $\lambda_0 = \lambda/(1-\lambda)$, then*

$$\lim_{n \to \infty} P\left\{ \frac{S_{n,a_n} + n^2 \ln(1 - a_n/n)}{\lambda_0^{1/2} n^{3/2}} \leq x \right\} = \Phi(x). \tag{3.3}$$

(c) *If $a_n/n \to 1$ and $b_n = n - a_n \to \infty$ as $n \to \infty$, then*

$$\lim_{n \to \infty} P\left\{ \frac{S_{n,a_n} + n^2 \ln(1 - a_n/n)}{n^2 b_n^{-1/2}} \leq x \right\} = \Phi(x). \tag{3.4}$$

**Theorem 3.3.** *If $n - a_n + 1 = k$ for a fixed positive integer $k$ and all positive integers $n$, then the limiting distribution of the random variable $S_{n,a_n} = M_n^{(k)}$ is given by the following equality:*

$$\lim_{n \to \infty} P\{M_n^{(k)} \leq n^2(x + \ln n)\} = \mathrm{e}^{-\mathrm{e}^{-x}} \sum_{s=0}^{k-1} \frac{\mathrm{e}^{-sx}}{s!}. \tag{3.5}$$

*In particular, the maximum $M_n = S_{nn}$ has, asymptotically, the Gumbel extreme value distribution.*

## 4. Proofs of Theorems 2.1, 2.2 and 2.3

**Proof of Theorem 2.1.** (a) It is easy to check that equality (2.10) holds for $k = 2$ and $k = 3$. The event $\{X_{nj} = k\}$, where $k > 3$, means that no two adjacent of the random variables $Z_1, Z_2, \ldots, Z_{k-3}$ take the value $j$ and that $Z_{k-2} \neq j$, $Z_{k-1} = Z_k = j$. Denote by $A_s$ the event that exactly $s$ of the random variables $Z_1, Z_2, \ldots, Z_{k-3}$ take the value $j$ and no two adjacent of them take the value $j$. Then

$$P(A_s) = \binom{k-s-2}{s}\left(1 - \frac{1}{n}\right)^{k-3-s} \frac{1}{n^s}, \qquad s \in \{0, 1, \ldots, [k/2] - 1\}. \tag{4.1}$$

The following two equalities hold:

$$\{X_{nj} = k\} = \bigcup_{s=0}^{[k/2]-1} \{A_s, Z_{k-2} \neq j, Z_{k-1} = Z_k = j\}, \tag{4.2}$$

$$P\{Z_{k-2} \neq j, Z_{k-1} = Z_k = j\} = \left(1 - \frac{1}{n}\right)\frac{1}{n^2}. \tag{4.3}$$

If $k > 3$, then the equality (2.10) follows from (4.1), (4.2) and (4.3).

(b) Using (2.10), we obtain that

$$1 - F_n(m) = \sum_{k=m+1}^{\infty} \frac{1}{n^2}\left(1 - \frac{1}{n}\right)^{k-2} \sum_{s=0}^{[(k-2)/2]} \binom{k-2-s}{s}\left(\frac{1}{n-1}\right)^s. \tag{4.4}$$

From problem 7(d), page 76 of Riordan (1968), the following identity holds:

$$\sum_{s=0}^{[r/2]} \binom{r-s}{s} x^s = \frac{1}{\alpha}\left\{\left(\frac{1+\alpha}{2}\right)^{r+1} - \left(\frac{1-\alpha}{2}\right)^{r+1}\right\}, \tag{4.5}$$

where $\alpha = (1 + 4x)^{1/2}$. The sum on the left-hand side of identity (4.5) is related to the Chebyshev polynomials. For $x = 1/(n-1)$, we get

$$\alpha = \left(1 + \frac{4}{n-1}\right)^{1/2} = \left(1 + \frac{3}{n}\right)^{1/2}\left(1 - \frac{1}{n}\right)^{-1/2} \tag{4.6}$$

and the tail $1 - F_n(m)$ can be represented in the form

$$1 - F_n(m) = \frac{1}{n^2}\left(1 - \frac{1}{n}\right)^{-1/2}\left(1 + \frac{3}{n}\right)^{-1/2}\left\{\frac{q_1^m}{1 - q_1} - \frac{q_2^m}{1 - q_2}\right\}, \tag{4.7}$$

where

$$q_1 = \left(1 - \frac{1}{n}\right)\frac{1+\alpha}{2} = 1 - \frac{1}{n^2} + \frac{5}{32n^3} + o\left(\frac{1}{n^3}\right), \tag{4.8}$$

$$q_2 = \left(1 - \frac{1}{n}\right)\frac{1-\alpha}{2} = -\frac{1}{n} + \frac{1}{n^2} - \frac{5}{32n^3} + o\left(\frac{1}{n^3}\right). \tag{4.9}$$

Let us determine $m$ from the condition $n(1 - F_n(m)) \to e^{-x}$ as $n \to \infty$. Using (4.7), (4.8) and (4.9), this condition can be transformed in the following way:

$$-\ln n - \frac{m}{n^2} + \frac{5m}{32n^3} + 2\ln n = -x + o(1) \quad \text{as } n \to \infty, \tag{4.10}$$

$$m = n^2(x + \ln n + o(1)) \quad \text{as } n \to \infty. \tag{4.11}$$

Consequently, (2.11) holds for $u_n = n^2(x + \ln n)$. □

**Proof of Theorem 2.2.** (a) Let $A$ be a subset of $N_n$, $|A| = j$ and let

$$a_0 = 1, \qquad b_0 = 0. \tag{4.12}$$

For any positive integer $l$, let us consider the set $S$ of all sequences of the form

$$c_1 c_2 \ldots c_l, \qquad \text{where } c_1, c_2, \ldots, c_l \in N_n, \tag{4.13}$$

such that no sequence from $S$ contains a subsequence of the form $aa$, $a \in A$. Let $a_l$ be the number of sequences from $S$ for which $c_l \in N_n \setminus A$, and $b_l$ the number of sequences from $S$ such that $c_l \in A$. The following equalities then hold:

$$a_1 = n - j = (n-j)(a_0 + b_0); \tag{4.14}$$

$$b_1 = j = ja_0 + (j-1)b_0; \tag{4.15}$$

$$a_{k-1} = (n-j)(a_{k-2} + b_{k-2}), \qquad \text{for any } k \geq 2; \tag{4.16}$$

$$b_{k-1} = ja_{k-2} + (j-1)b_{k-2}, \qquad \text{for any } k \geq 2. \tag{4.17}$$

Let $s_l = a_l + b_l$ for $l \geq 0$. It follows from (4.12) and (4.14)–(4.17) that

$$s_0 = a_0 + b_0 = 1, \qquad s_1 = a_1 + b_1 = n, \tag{4.18}$$

$$s_{k-1} = a_{k-1} + b_{k-1} = (n-1)s_{k-2} + (n-j)s_{k-3} \qquad \text{for } k \geq 3. \tag{4.19}$$

Hence, the sequence $(s_l)$ satisfies the linear difference equation (4.19) with initial conditions (4.18). It follows that

$$s_{k-1} = C_1 t_1^{k-1} + C_2 t_2^{k-1} \qquad \text{for any } k \geq 1, \tag{4.20}$$

where $t_1 = t_1(j)$ and $t_2 = t_2(j)$ are given by (2.13) and (2.14). Using initial conditions (4.18), we obtain the constants $C_1$ and $C_2$:

$$C_1 = \frac{1}{2}\left\{1 + \left(1 - \frac{4j}{(n+1)^2}\right)^{-1/2}\right\}, \tag{4.21}$$

$$C_2 = \frac{1}{2}\left\{1 - \left(1 - \frac{4j}{(n+1)^2}\right)^{-1/2}\right\}. \tag{4.22}$$

Note that

$$P\{Y_{nj} = k\} = b_{k-1} n^{-k}, \qquad k \in \{2, 3, \ldots\}, \tag{4.23}$$

$$b_{k-1} = s_{k-1} - a_{k-1} = s_{k-1} - (n-j)s_{k-2}. \tag{4.24}$$

Equalities (2.12) follow from (4.20)–(4.24).

(b) Let $D = \{(n+1)^2 - 4j\}^{1/2}$. We then have

$$\mathrm{E}Y_{nj} = \frac{j}{nD} \cdot \left\{ \sum_{k=2}^{\infty} k \left( \frac{t_1}{n} \right)^{k-1} - \sum_{k=2}^{\infty} k \left( \frac{t_2}{n} \right)^{k-1} \right\}. \tag{4.25}$$

Note that $\sum_{k=2}^{\infty} k q^{k-1} = \frac{2q-q^2}{(1-q)^2}$ for $|q| < 1$. We now get

$$\mathrm{E}Y_{nj} = \frac{j}{nD} \cdot \left\{ \frac{t_1(2n-t_1)}{(n-t_1)^2} - \frac{t_2(2n-t_2)}{(n-t_2)^2} \right\}. \tag{4.26}$$

Since $t_1 = (n-1+D)/2$, $t_2 = (n-1-D)/2$, $2n-t_1 = (3n+1-D)/2$, $2n-t_2 = (3n+1+D)/2$, $n-t_1 = (n+1-D)/2$ and $n-t_2 = (n+1+D)/2$, the mean $\mathrm{E}Y_{nj}$ can be represented in the form

$$\mathrm{E}Y_{nj} = \frac{j}{nD} \cdot \frac{U_1}{V_1}, \tag{4.27}$$

where $U_1$ and $U_2$ can be transformed in the following way:

$$\begin{aligned}
U_1 &= (n-1+D)(n+1+D)^2(3n+1-D) \\
&\quad - (n-1-D)(n+1-D)^2(3n+1+D) \\
&= (n^2 + 2nD + D^2 - 1)(3n^2 + 4n + 2nD + 1 - D^2) \\
&\quad - (n^2 - 2nD + D^2 - 1)(3n^2 + 4n - 2nD + 1 - D^2) \\
&= 16nD(n^2+n).
\end{aligned}$$

It also follows that $V_1 = (n+1-D)^2(n+1+D)^2 = 16j^2$. The first of the equalities (2.15) follows easily from (4.27). Let us now determine $\mathrm{var}\, Y_{nj}$. We have

$$\mathrm{E}(Y_{nj}^2) = \frac{j}{nD} \cdot \left\{ \sum_{k=2}^{\infty} k^2 \left( \frac{t_1}{n} \right)^{k-1} - \sum_{k=2}^{\infty} k^2 \left( \frac{t_2}{n} \right)^{k-1} \right\}. \tag{4.28}$$

Since $\sum_{k=2}^{\infty} k^2 q^{k-1} = \frac{q^3 - 3q^2 + 4q}{(1-q)^3}$ for $|q| < 1$, we obtain

$$\begin{aligned}
\mathrm{E}Y_{nj}^2 &= \frac{j}{nD} \cdot \left\{ \frac{t_1^3 - 3nt_1^2 + 4n^2 t_1}{(n-t_1)^2} - \frac{t_1^3 - 3nt_1^2 + 4n^2 t_1}{(n-t_1)^2} \right\} \\
&= \frac{j}{nD} \cdot \left\{ \frac{(n-1+D)^3 - 6(n-1+D)^2 + 16n^2(n-1+D)}{(n+1-D)^3} \right. \\
&\qquad\qquad \left. - \frac{(n-1-D)^3 - 6n(n-1-D)^2 + 16n^2(n-1-D)}{(n+1+D)^3} \right\} = \frac{j}{nD} \cdot \frac{U_2}{V_2},
\end{aligned}$$

where $V_2 = (n+1-D)^3(n+1+D)^3 = 64j^3$ and

$$U_2 = \{(n-1+D)^3 - 6n(n-1+D)^2 + 16n^2(n-1+D)\} \cdot (n+1+D)^3$$

$$- \{(n-1-D)^3 - 6n(n-1-D)^2 + 16n^2(n-1-D)\} \cdot (n+1-D)^3$$
$$= (n^2 + 2nD + D^2 - 1)^3 - (n^2 - 2nD + D^2 - 1)^3$$
$$- 6n\{(n+1+D)(n^2+2nD+D^2-1)^2$$
$$- (n+1-D)(n^2-2nD+D^2-1)^2\}$$
$$+ 16n^2\{(n+1+D)^2(n^2+2nD+D^2-1)$$
$$- (n+1-D)^2(n^2-2nD+D^2-1)\}.$$

Let $M = n^2 + D^2 - 1 = 2n^2 + 2n - 4j$. Since $D^2 = (n+1)^2 - 4j$, we get

$$U_2 = 4nD\{(M+2nD)^2 + (M+2nD)(M-2nD) + (M-2nD)^2\}$$
$$- 6n\{(n+1+D)(M+2nD)^2 - (n+1-D)(M-2nD)^2\}$$
$$+ 16n^2\{(n+1+D)^2(M+2nD) - (n+1-D)^2(M-2nD)\}$$
$$= 4nD(3M^2 + 4n^2D^2) - 12nD\{4n(n+1)M + M^2 + 4n^2D^2\}$$
$$+ 64n^2D\{(n+1)M + n(n+1)^2 + nD^2\}$$
$$= 64nD(2n^4 + 4n^3 + 2n^2 - 3n^2j - nj)$$

and, consequently,

$$\mathrm{E} Y_{nj}^2 = \frac{j}{nD} \cdot \frac{U_2}{V_2} = \frac{1}{j^2}(2n^4 + 4n^3 + 2n^2 - 3n^2 j - nj). \tag{4.29}$$

The second of the equalities (2.15) follows from (4.29) and $\mathrm{E} Y_{nj} = (n^2+n)/j$. □

**Proof of Theorem 2.3.** For any positive integer $n$, let

$$H_n = 1 + \frac{1}{2} + \cdots + \frac{1}{n}, \qquad H_n^{(2)} = 1 + \frac{1}{2^2} + \cdots + \frac{1}{n^2}. \tag{4.30}$$

The equalities

$$H_n = \ln n + \gamma + \frac{1}{2n} - \frac{1}{12n^2} + \frac{1}{120n^4} - \frac{\varepsilon_n}{252n^6}, \tag{4.31}$$

$$H_n^{(2)} = \frac{\pi^2}{6} - \frac{1}{n} + \frac{\vartheta_n}{n(n+1)}, \tag{4.32}$$

hold, where $0 < \varepsilon_n < 1$ and $0 < \vartheta_n < 1$. The equality (4.31) can be found in Graham, Knuth and Patashnik (1994), page 480, and (4.32) can easily be proven. Since

$$\mathrm{E} S_{n,a_n} = (n^2 + n)(H_n - H_{n-a_n}), \tag{4.33}$$

relations concerning the asymptotic behavior of $\mathrm{E} S_{n,a_n}$ follow from (4.31) and (4.33). Similarly, using (4.31), (4.32) and the second of the equalities (2.15), we get relations concerning the asymptotic behavior of $\mathrm{var}\, S_{n,a_n}$. □

## 5. Proofs of Theorems 3.1, 3.2 and 3.3

Let $D = D(j) = \{(n+1)^2 - 4j\}^{1/2}$ and let $t_1 = t_1(j)$ and $t_2 = t_2(j)$ be given by (2.13) and (2.14). In the sequel, we shall use the following series expansions and approximations that follow from them:

$$D = n\left(1 + \frac{1}{n} - \frac{2j}{n^2} + \frac{2j}{n^3} - \frac{2j+2j^2}{n^4} + \frac{2j+6j^2}{n^5} - \cdots\right), \tag{5.1}$$

$$D^{-1} = \frac{1}{n}\left(1 - \frac{1}{n} + \frac{1+2j}{n^2} - \frac{1+6j}{n^3} + \cdots\right), \tag{5.2}$$

$$t_1 = n\left(1 - \frac{j}{n^2} + \frac{j}{n^3} - \frac{j+j^2}{n^4} + \frac{j+3j^2}{n^5} - \cdots\right), \tag{5.3}$$

$$t_2 = -1 + \frac{j}{n} - \frac{j}{n^2} + \frac{j+j^2}{n^3} - \frac{j+3j^2}{n^4} - \cdots, \tag{5.4}$$

$$n - t_1 = \frac{j}{n}\left(1 - \frac{1}{n} + \frac{1+j}{n^2} - \frac{1+3j}{n^3} - \cdots\right). \tag{5.5}$$

**Proof of Theorem 3.1.** Let us determine the characteristic function of the random variable $Y_{nj}$. Using (2.12), we obtain

$$f_{nj}(t) = \frac{j}{nD}\left\{\sum_{k=2}^{\infty} e^{itk}\left(\frac{t_1}{n}\right)^{k-1} - \sum_{k=2}^{\infty} e^{itk}\left(\frac{t_2}{n}\right)^{k-1}\right\}$$

$$= \frac{j}{t_1 D} \cdot \sum_{k=2}^{\infty}\left(e^{it}\frac{t_1}{n}\right)^k - \frac{j}{t_2 D} \cdot \sum_{k=2}^{\infty}\left(e^{it}\frac{t_2}{n}\right)^k$$

$$= \frac{j}{t_1 D} \cdot \frac{e^{2it}t_1^2/n^2}{1 - e^{it}t_1/n} - \frac{j}{t_2 D} \cdot \frac{e^{2it}t_2^2/n^2}{1 - e^{it}t_2/n}$$

$$= \frac{j}{t_1 D} \cdot \frac{e^{2it}t_1^2/n^2}{1 - e^{it}t_1/n}\left\{1 - \frac{t_2}{t_1} \cdot \frac{1 - e^{it}t_1/n}{1 - e^{it}t_2/n}\right\}.$$

Consequently, we get that the characteristic function of the random variable $Y_{nj}$ can be represented in the form

$$f_{nj}(t) = \frac{j}{t_1}\left(\frac{t_1}{n}\right)^2\left(1 + \frac{2}{n} + \frac{1-4j}{n^2}\right)^{-1/2}\left\{1 - \frac{t_2}{t_1} \cdot \frac{1-e^{it}t_1/n}{1-e^{it}t_2/n}\right\}\frac{e^{2it}}{n - t_1 e^{it}}. \tag{5.6}$$

If $j/n \to 1$ as $n \to \infty$, then

$$\lim_{n \to \infty} \frac{j}{t_1}\left(\frac{t_1}{n}\right)^2\left(1 + \frac{2}{n} + \frac{1-4j}{n^2}\right)^{-1/2}\left\{1 - \frac{t_2}{t_1} \cdot \frac{1-e^{it/n}t_1/n}{1-e^{it/n}t_2/n}\right\} = 1 \tag{5.7}$$

and hence the asymptotic behavior of $f_{Y_{nj}/n}(t)$ is given by

$$f_{Y_{nj}/n}(t) \sim \frac{e^{2it/n}}{n - t_1 e^{it/n}}, \qquad n \to \infty. \tag{5.8}$$

The relations

$$|n - t_1 e^{it/n}|^2 = \left|n - t_1 \cos\frac{t}{n} - it_1 \sin\frac{t}{n}\right|^2 = n^2 + t_1^2 - 2nt_1 \cos\frac{t}{n}$$

$$= (n - t_1)^2 + 2nt_1\left(1 - \cos\frac{t}{n}\right) \to 1 + t^2 \qquad \text{as } n \to \infty, j/n \to 1,$$

$$\arg f_{Y_{nj}/n}(t) \sim \arg\frac{e^{2it}}{n - t_1 e^{it/n}} = \frac{2t}{n} - \arctan\frac{-t_1 \sin(t/n)}{n - t_1 \cos(t/n)}$$

$$= \frac{2t}{n} - \arctan\frac{-t + o(1)}{1 + o(1)} \to \arctan t \qquad \text{as } n \to \infty$$

hold, where we again used the assumption that $j/n \to 1$ as $n \to \infty$. Hence, $f_{Y_{nj}/n}(t) \to (1+t^2)^{-1/2} e^{i \cdot \arctan t}$ as $n \to \infty$. Consequently, if $a_n = k$, where $k$ is a fixed positive integer, we get that $(Y_{nn} + \cdots + Y_{n,n-k+1})/n$ converges in distribution to a random variable whose characteristic function is given by (3.1). $\square$

**Proof of Theorem 3.2.** Using (5.6), we obtain that the characteristic function of the sum $S_{n,a_n}$ can be represented in the following way:

$$f_{S_{n,a_n}}(t) = \prod_{j=n-a_n+1}^{n} f_{nj}(t) = P_{n1} \cdot P_{n2} \cdot P_{n3}(t) \cdot P_{n4}(t), \tag{5.9}$$

where $P_{n1}$, $P_{n2}$, $P_{n3}(t)$ and $P_{n4}(t)$ are given by

$$P_{n1} = \prod_{j=n-a_n+1}^{n} \frac{t_1(j)}{n}, \tag{5.10}$$

$$P_{n2} = \prod_{j=n-a_n+1}^{n} \left(1 + \frac{2}{n} + \frac{1-4j}{n^2}\right)^{-1/2}, \tag{5.11}$$

$$P_{n3}(t) = \prod_{j=n-a_n+1}^{n} \left\{1 - \frac{t_2(j)}{t_1(j)} \cdot \frac{1 - e^{it} \cdot t_1(j)/n}{1 - e^{it} \cdot t_2(j)/n}\right\}, \tag{5.12}$$

$$P_{n4}(t) = \prod_{j=n-a_n+1}^{n} \frac{e^{2it}}{\frac{n^2}{j}(1 - (t_1/n)e^{it})}. \tag{5.13}$$

**Lemma 5.1.** *If $a_n/n \to \lambda \in [0,1]$ as $n \to \infty$, then the following relation holds:*

$$\lim_{n\to\infty} P_{n1} = e^{-\lambda + \lambda^2/2}. \tag{5.14}$$

**Proof.** The relations

$$1 - \frac{j}{n^2} + \frac{c_1}{n^2} \le \frac{t_1(j)}{n} \le 1 - \frac{j}{n^2} + \frac{c_2}{n^2}, \tag{5.15}$$

$$\ln(1+x) = x + r_1(x), \qquad |r_1(x)| \le |x|^2, \qquad \text{for } |x| < 1, \tag{5.16}$$

hold, where constants $c_1$ and $c_2$ do not depend on $j$. Consequently, we get

$$\ln P_{n1} \sim -\frac{1}{n^2} \sum_{j=n-a_n+1}^{n} j = -\frac{1}{n^2} \cdot \frac{(2n - a_n + 1)a_n}{2} = -\frac{a_n}{n} + \frac{a_n^2}{2n^2} - \frac{a_n}{2n^2} \tag{5.17}$$

and the statement of the lemma follows easily. $\square$

**Lemma 5.2.** *If $a_n/n \to \lambda \in [0,1]$ as $n \to \infty$, then the following relation holds:*

$$\lim_{n\to\infty} P_{n2} = e^{\lambda - \lambda^2}. \tag{5.18}$$

**Proof.** The statement of the lemma follows from the following relations:

$$\ln P_{n2} = -\frac{1}{2} \sum_{j=n-a_n+1}^{n} \ln\left(1 + \frac{2}{n} + \frac{1 - 4j}{n^2}\right) \sim -\frac{1}{2} \sum_{j=n-a_n+1}^{n} \frac{2n - 4j}{n^2}$$

$$= -\frac{1}{n^2} \sum_{j=n-a_n+1}^{n} (n - 2j) = \frac{a_n}{n} - \frac{a_n^2}{n^2} + \frac{a_n}{n^2}.$$

$\square$

**Lemma 5.3.** *If $\tau_n^2$ is the main term of $\sigma_n^2$ determined by (2.16)–(2.18), $\tau_n > 0$ and $a_n \to \infty$ as $n \to \infty$, then for any real $t$ the following equality holds:*

$$\lim_{n\to\infty} P_{n3}\left(\frac{t}{\tau_n}\right) = 1. \tag{5.19}$$

**Proof.** Using (2.13) and (2.14), it is easy to prove that the following inequalities hold for any $j \in \{1, 2, \ldots, n\}$:

$$-\frac{1}{n} \le \frac{t_2(j)}{t_1(j)} \le 0 \le 1 - \frac{t_1(j)}{n} \le \frac{1}{n}. \tag{5.20}$$

For sufficiently large $n$, the inequality $1 - \frac{1}{n} \le \cos\frac{t}{\tau_n} \le 1$ also holds and for such values of $n$, we obtain

$$\frac{|1 - \mathrm{e}^{\mathrm{i}t/\tau_n} \cdot (t_1/n)|^2}{|1 - \mathrm{e}^{\mathrm{i}t/\tau_n} \cdot (t_2/n)|^2} = \frac{1 - (2t_1/n)\cos(t/\tau_n) + t_1^2/n^2}{1 - (2t_2/n)\cos(t/\tau_n) + t_2^2/n^2}$$

$$\leq \frac{1 + (2t_1/n)(1/n - 1) + t_1^2/n^2}{1 - 2t_2/n + t_2^2/n^2} = \frac{(1 - t_1/n)^2 + 2t_1/n^2}{(1 - t_2/n)^2} \quad (5.21)$$

$$\leq \frac{1}{n^2} + \frac{2}{n} \leq \frac{3}{n}.$$

Using the first of the inequalities (5.19) and the inequality (5.20), we obtain from (5.12) that $\arg P_{n3}(t/\tau_n) \to 0$ and $P_{n3}(t/\tau_n) \to 1$ as $n \to \infty$ for every real $t$. □

**Lemma 5.4.** *Let $\tau_n^2$ be the main term of $\sigma_n^2$, $\tau_n > 0$ and suppose that $a_n \to \infty$ and $b_n = n - a_n \to \infty$ as $n \to \infty$. For any real $t$, the following asymptotic relations then hold as $n \to \infty$:*

$$\ln \left| P_{n4}\left(\frac{t}{\tau_n}\right) \right| \sim \frac{a_n^2}{2n^2} - \frac{n^4 t^2}{2\tau_n^2}(H_n^{(2)} - H_{n-a_n}^{(2)}); \quad (5.22)$$

$$\left| P_{n4}\left(\frac{t}{\tau_n}\right) \right| \to \exp\left(\frac{\lambda^2}{2} - \frac{t^2}{2}\right) \quad \text{if } \frac{a_n}{n} \to \lambda \in [0, 1]. \quad (5.23)$$

**Proof.** (a) We shall use the following inequalities:

$$1 - \frac{1}{n} \leq \frac{t_1}{n} \leq 1, \quad (5.24)$$

$$\frac{t^2}{2\tau_n^2} - \frac{t^4}{24\tau_n^4} \leq 1 - \cos\frac{t}{\tau_n} \leq \frac{t^2}{2\tau_n^2}, \quad (5.25)$$

$$1 - \frac{2}{n} + \frac{c_1}{n^2} \leq \frac{n^4}{j^2}\left(1 - \frac{t_1}{n}\right)^2 \leq 1 - \frac{2}{n} + \frac{2j}{n^2} + \frac{c_2}{n^2}, \quad (5.26)$$

where the constants $c_1$ and $c_2$ do not depend on $j$. Inequalities (5.24) and (5.25) are straightforward exercises and (5.26) follows from the equality

$$\frac{n^4}{j^2}\left(1 - \frac{t_1}{n}\right)^2 = \frac{n^4}{2j^2 n^2}\left\{2(n+1)^2 - 4j - 2n(n+1)\left(1 + \frac{2}{n} + \frac{1-4j}{n^2}\right)^{1/2}\right\}. \quad (5.27)$$

Using the equality

$$\left| 1 - \frac{t_1}{n}\mathrm{e}^{\mathrm{i}t/\tau_n} \right|^2 = \left(1 - \frac{t_1}{n}\right)^2 + \frac{2t_1}{n}\left(1 - \cos\frac{t}{\tau_n}\right) \quad (5.28)$$

and inequalities (5.24)–(5.26), we get

$$1 - \frac{2}{n} + \frac{2j}{n^2} + \frac{c_1}{n^2} + \frac{2n^4}{j^2}\left(1 - \frac{1}{n}\right)\left(\frac{t^2}{2\tau_n^2} - \frac{t^4}{24\tau_n^4}\right) \leq \frac{n^4}{j^2}\left| 1 - \frac{t_1}{n}\mathrm{e}^{\mathrm{i}t/\tau_n} \right|^2 \quad (5.29)$$

$$\frac{n^4}{j^2}\left|1 - \frac{t_1}{n}e^{it/\tau_n}\right|^2 \leq 1 - \frac{2}{n} + \frac{2j}{n^2} + \frac{c_2}{n^2} + \frac{n^4 t^2}{\tau_n^2 j^2}. \quad (5.30)$$

It follows from (5.13) that

$$\left|P_{n4}\left(\frac{t}{\tau_n}\right)\right| = \left\{\prod_{j=n-a_n+1}^{n} \frac{n^4}{j^2}\left|1 - \frac{t_1(j)}{n}e^{it/\tau_n}\right|^2\right\}^{-1/2}. \quad (5.31)$$

Finally, using (5.16) and (5.29)–(5.31), we get, as $n \to \infty$,

$$\ln\left|P_{n4}\left(\frac{t}{\tau_n}\right)\right| \sim \sum_{j=n-a_n+1}^{n} \ln\left(1 - \frac{2}{n} + \frac{2j}{n^2} + \frac{n^4 t^2}{\tau_n^2 j^2}\right)^{-1/2}$$

$$\sim \frac{1}{2} \sum_{j=n-a_n+1}^{n} \left(\frac{2}{n} - \frac{2j}{n^2} - \frac{n^4 t^2}{\tau_n^2 j^2}\right)$$

$$\sim \frac{a_n^2}{2n^2} - \frac{n^4 t^2}{2\tau_n^2}(H_n^{(2)} - H_{n-a_n}^{(2)}).$$

(b) Using (4.32) and the main term $\tau_n^2$ of variance $\sigma_n^2$ from relations (2.16)–(2.18), we get relation (5.23) □

**Lemma 5.5.** *If $\nu_n$ and $\tau_n^2$ are the main terms of the mean $\mu_n = \mathrm{E} S_{n,a_n}$ and the variance $\sigma_n^2 = \mathrm{var}\, S_{n,a_n}$, $\tau_n > 0$ and $S_{n,a_n}^* = \tau_n^{-1}(S_{n,a_n} - \nu_n)$, then*

$$\arg f_{S_{n,a_n}^*}(t) = o(1) \qquad \text{as } n \to \infty. \quad (5.32)$$

**Proof.** The following equalities hold:

$$\arg P_{n4}\left(\frac{t}{\tau_n}\right) = \frac{2ta_n}{\tau_n} - \sum_{j=n-a_n+1}^{n} \arg\left(1 - \frac{t_1(j)}{n}e^{it/\tau_n}\right)$$

$$= \frac{2ta_n}{\tau_n} + \sum_{j=n-a_n+1}^{n} \arctan\frac{N_n(t,j)}{D_n(t,j)}, \quad (5.33)$$

where

$$N_n(t,j) = \frac{t_1(j)}{n}\sin\frac{t}{\tau_n}, \qquad D_n(t,j) = 1 - \frac{t_1(j)}{n}\cos\frac{t}{\tau_n}. \quad (5.34)$$

*Case* 1. Let $a_n \to \infty$, $\frac{a_n}{n} \to 0$ as $n \to \infty$. We then have $\tau_n^2 = n^2 a_n$, $\tau_n = na_n^{1/2}$ and $\nu_n = -n^2 \ln(1 - \frac{a_n}{n})$. For $N_n(t,j)$ and $D_n^{-1}(t,j)$, we obtain

$$N_n(t,j) = \frac{t}{na_n^{1/2}}\left\{1 - \frac{j}{n^2} + \frac{\vartheta C}{n^2}\right\}, \quad (5.35)$$

$$D_n^{-1}(t,j) = \frac{n^2}{j}\left\{1 + \frac{1}{n} - \frac{j}{n^2} - \frac{t^2}{2ja_n} + \frac{\vartheta C}{n^2}\right\}, \tag{5.36}$$

$$\frac{N_n(t,j)}{D_n(t,j)} = \frac{tn}{ja_n^{1/2}}\left\{1 + \frac{1}{n} - \frac{2j}{n^2} - \frac{t^2}{2ja_n} + \frac{\vartheta C}{n^2}\right\}. \tag{5.37}$$

In (5.35)–(5.37) and in relations that will follow $C = C(t) > 0$ is a constant which does not depend on $j$, and $\vartheta \in [-1,1]$. Also, note that $\vartheta$ may be different at different occurrences. Consequently, we obtain the following results:

$$\arg P_{n4}\left(\frac{t}{\tau_n}\right) = \frac{2ta_n}{na_n^{1/2}} + \sum_{j=n-a_n+1}^{n} \frac{tn}{ja_n^{1/2}}\left\{1 + \frac{1}{n} - \frac{2j}{n^2} - \frac{t^2}{2ja_n} + \frac{\vartheta C}{n^2}\right\}$$

$$= \frac{2ta_n^{1/2}}{n} + \frac{tn}{a_n^{1/2}}\left(1 + \frac{1}{n}\right)(H_n - H_{n-a_n})$$

$$- \frac{2ta_n^{1/2}}{n} - \frac{t^3 n}{2a_n^{3/2}}(H_n^{(2)} - H_{n-a_n}^{(2)}) + o(1)$$

$$= -\frac{tn}{a_n^{1/2}}\left(1 + \frac{1}{n}\right)\ln\left(1 - \frac{a_n}{n}\right) - \frac{t^3 n}{2a_n^{3/2}} \cdot \frac{a_n}{n^2}\left(1 - \frac{a_n}{n}\right)^{-1} + o(1)$$

$$= ta_n^{1/2} - \frac{t^3}{2na_n^{1/2}}\left(1 - \frac{a_n}{n}\right)^{-1} + o(1)$$

$$= ta_n^{1/2} + o(1) \quad \text{as } n \to \infty;$$

$$\arg \prod_{j=n-a_n+1}^{n} f_{nj}\left(\frac{t}{na_n^{1/2}}\right) = ta_n^{1/2} + o(1), \quad n \to \infty;$$

$$f_{S_{n,a_n}^*}(t) = \exp\left(-\frac{it\nu_n}{na_n^{1/2}}\right)\prod_{j=n-a_n+1}^{n} f_{nj}\left(\frac{t}{na_n^{1/2}}\right);$$

$$\arg f_{S_{n,a_n}^*}(t) = -\frac{t\nu_n}{na_n^{1/2}} + ta_n^{1/2} + o(1) = o(1), \quad n \to \infty. \tag{5.38}$$

*Case* 2. Let $\frac{a_n}{n} \to \lambda \in (0,1)$ as $n \to \infty$, and $\lambda_0 = \frac{\lambda}{1-\lambda}$. We then have $\tau_n^2 = \lambda_0 n^3$, $\tau_n = \lambda_0^{1/2} n^{3/2}$ and $\nu_n = -n^2 \ln(1 - \frac{a_n}{n})$. We now get

$$N_n(t,j) = \frac{t}{\lambda_0^{1/2} n^{3/2}}\left\{1 - \frac{j}{n^2} + \frac{\vartheta C}{n^2}\right\}, \tag{5.39}$$

$$D_n^{-1}(t,j) = \frac{n^2}{j}\left\{1 + \frac{1}{n} - \frac{1+j}{n^2} - \frac{t^2}{2\lambda_0 jn} + \frac{\vartheta C}{n^2}\right\}, \tag{5.40}$$

$$\frac{N_n(t,j)}{D_n(t,j)} = \frac{tn^{1/2}}{\lambda_0^{1/2} j}\left\{1 + \frac{1}{n} - \frac{2j}{n^2} - \frac{t^2}{2\lambda_0 jn} + \frac{\vartheta C}{n^2}\right\} \tag{5.41}$$

and, consequently, we obtain the following:

$$\arg P_{n4}\left(\frac{t}{\tau_n}\right) = \frac{2ta_n}{\lambda_0^{1/2} n^{3/2}} + \sum_{j=n-a_n+1}^{n} \frac{tn^{1/2}}{\lambda_0^{1/2} j}\left\{1 + \frac{1}{n} - \frac{2j}{n^2} - \frac{t^2}{2\lambda_0 jn} + \frac{\vartheta C}{n^2}\right\}$$

$$= \frac{2ta_n}{\lambda_0^{1/2} n^{3/2}} + \frac{tn^{1/2}}{\lambda_0^{1/2}}\left(1 + \frac{1}{n}\right)(H_n - H_{n-a_n})$$

$$- \frac{2ta_n}{\lambda_0^{1/2} n^{3/2}} - \frac{t^3}{2\lambda_0^{3/2} n^{1/2}}(H_n^{(2)} - H_{n-a_n}^{(2)}) + o(1)$$

$$= \frac{tn^{1/2}}{\lambda_0^{1/2}}\left(1 + \frac{1}{n}\right)(H_n - H_{n-a_n}) + o(1)$$

$$= \frac{tn^{1/2}}{\lambda_0^{1/2}}\left(1 + \frac{1}{n}\right)$$

$$\times \left\{\ln n - \ln(n-a_n) + \frac{1}{2n} - \frac{1}{2(n-a_n)} + \frac{\vartheta C}{n^2}\right\} + o(1)$$

$$= -\frac{tn^{1/2}}{\lambda_0^{1/2}}\ln\left(1 - \frac{a_n}{n}\right) + o(1) \qquad \text{as } n \to \infty; \tag{5.42}$$

$$\arg \prod_{j=n-a_n+1}^{n} f_{nj}\left(\frac{t}{\lambda_0^{1/2} n^{3/2}}\right) = -\frac{tn^{1/2}}{\lambda_0^{1/2}}\ln\left(1 - \frac{a_n}{n}\right) + o(1), \qquad n \to \infty;$$

$$f_{S_{n,a_n}^*}(t) = \exp\left(-\frac{it\nu_n}{\lambda_0^{1/2} n^{3/2}}\right) \prod_{j=n-a_n+1}^{n} f_{nj}\left(\frac{t}{\lambda_0^{1/2} n^{3/2}}\right);$$

$$\arg f_{S_{n,a_n}^*}(t) = -\frac{t\nu_n}{\lambda_0^{1/2} n^{3/2}} - \frac{tn^{1/2}}{\lambda_0^{1/2}}\ln\left(1 - \frac{a_n}{n}\right) + o(1) = o(1), \qquad n \to \infty.$$

*Case* 3. Let $\frac{a_n}{n} \to 1$ as $n \to \infty$ and $b_n = n - a_n > 0$ for all $n$. We then have $\tau_n^2 = n^4 b_n^{-1}$, $\tau_n = n^2 b_n^{-1/2}$ and $\nu_n = -n^2 \ln(1 - \frac{a_n}{n}) = n^2 \ln \frac{n}{b_n}$. We now get

$$N_n(t,j) = \frac{tb_n^{1/2}}{n^2}\left\{1 - \frac{j}{n^2} + \frac{\vartheta C}{n^2}\right\}, \tag{5.43}$$

$$D_n^{-1}(t,j) = \frac{n^2}{j}\left\{1 + \frac{1}{n} - \frac{j}{n^2} - \frac{t^2 b_n}{2jn^2} + \frac{\vartheta C}{n^2}\right\}, \tag{5.44}$$

$$\frac{N_n(t,j)}{D_n(t,j)} = \frac{tb_n^{1/2}}{j}\left\{1 + \frac{1}{n} - \frac{2j}{n^2} - \frac{t^2 b_n}{2jn^2} + \frac{\vartheta C}{n^2}\right\} \tag{5.45}$$

and, consequently, we obtain the following:

$$\arg P_{n4}\left(\frac{t}{\tau_n}\right) = \frac{2ta_n b_n^{1/2}}{n^2} + \sum_{j=n-a_n+1}^{n} \frac{tb_n^{1/2}}{j}\left\{1 + \frac{1}{n} - \frac{2j}{n^2} - \frac{t^2 b_n}{2jn^2} + \frac{\vartheta C}{n^2}\right\}$$

$$= \frac{2ta_n b_n^{1/2}}{n^2} + tb_n^{1/2}\left(1 + \frac{1}{n}\right)(H_n - H_{n-a_n})$$

$$- \frac{2ta_n b_n^{1/2}}{n^2} - \frac{t^3 b_n^{3/2}}{2n^2}(H_n^{(2)} - H_{n-a_n}^{(2)}) + o(1)$$

$$= tb_n^{1/2}\left(1 + \frac{1}{n}\right)(H_n - H_{n-a_n}) + o(1)$$

$$= tb_n^{1/2}\left(1 + \frac{1}{n}\right)\left(\ln n - \ln b_n + \frac{1}{2n} - \frac{1}{2b_n} + \cdots\right) + o(1)$$

$$= tb_n^{1/2} \ln \frac{n}{b_n} + o(1) \quad \text{as } n \to \infty; \tag{5.46}$$

$$\arg \prod_{j=n-a_n+1}^{n} f_{nj}\left(\frac{tb_n^{1/2}}{n^2}\right) = tb_n^{1/2} \ln \frac{n}{b_n} + o(1), \quad n \to \infty;$$

$$f_{S_{n,a_n}^*}(t) = \exp\left(-\frac{it\nu_n}{n^2 b_n^{-1/2}}\right) \prod_{j=n-a_n+1}^{n} f_{nj}\left(\frac{tb_n^{1/2}}{n^2}\right);$$

$$\arg f_{S_{n,a_n}^*}(t) = -\frac{t\nu_n}{n^2 b_n^{-1/2}} + tb_n^{1/2} \ln \frac{n}{b_n} + o(1) = o(1), \quad n \to \infty.$$

Finally, using relations (5.9)–(5.14), (5.18), (5.19), (5.23) and (5.32), we obtain that in all cases considered, $f_{S_{n,a_n}^*}(t) \to e^{-t^2/2}$ as $n \to \infty$ and therefore the proof of Theorem 3.2 is completed. $\square$

**Proof of Theorem 3.3.** We shall first prove a few lemmas.

**Lemma 5.6.** *Let $X_{nj}^*$, $j \in \{1, 2, \ldots, n\}$, be independent random variables with the same probability distribution,*

$$P\{X_{nj}^* = k\} = \sum_{s=0}^{[k/2]-1} \binom{k-2-s}{s}\left(1 - \frac{1}{n}\right)^{k-s-2} \frac{1}{n^{s+2}}, \quad k \geq 2, \tag{5.47}$$

and let $M_n^* = \max\{X_{n1}^*, \ldots, X_{nn}^*\}$. For every real $x$, the following equality then holds:

$$\lim_{n \to \infty} P\{M_n^* \leq n^2(x + \ln n)\} = e^{-e^{-x}}. \tag{5.48}$$

**Proof.** The statements $P\{M_n^* \leq u_n\} \to e^{-\tau}$ as $n \to \infty$ and $n(1 - F_n(u_n)) \to \tau$ as $n \to \infty$ are equivalent. Hence, (5.48) follows from (2.11). □

**Lemma 5.7.** *If $j$ is a fixed positive integer and $n \to \infty$, then the following asymptotic relation holds:*

$$P\{Y_{nj} > n^2(x + \ln n)\} = \frac{e^{-jx}}{n^j}\left\{1 + \frac{j(x + \ln n)}{n} + O\left(\left(\frac{\ln n}{n}\right)^2\right)\right\}. \tag{5.49}$$

**Proof.** Let $u_n = n^2(x + \ln n)$ and $r_n = u_n - [u_n]$. Using (2.12), we obtain that

$$P\{Y_{nj} > u_n\} = P\{Y_{nj} > [u_n]\}$$

$$= \frac{j}{(n - t_1)\{(n+1)^2 - 4j\}^{1/2}} \cdot \left(\frac{t_1}{n}\right)^{[u_n]}$$

$$- \frac{j}{(n - t_2)\{(n+1)^2 - 4j\}^{1/2}} \cdot \left(\frac{t_2}{n}\right)^{[u_n]}$$

$$\equiv A_1 - A_2.$$

Using (5.1), (5.3) and (5.5), we obtain that

$$\ln A_1 = \ln j - \ln(n - t_1) - \ln\{(n+1)^2 - 4j\}^{1/2} + (u_n - r_n)\ln(t_1/n)$$

$$= \ln j - \ln j + \ln n - \ln\left(1 - \frac{1}{n} + O\left(\frac{1}{n^2}\right)\right) - \ln n - \ln\left(1 + \frac{1}{n} + O\left(\frac{1}{n^2}\right)\right)$$

$$+ \{n^2(x + \ln n) - r_n\}\ln\left(1 - \frac{j}{n^2} + \frac{j}{n^3} + O\left(\frac{1}{n^4}\right)\right)$$

$$= -j(x + \ln n) + \frac{j(x + \ln n)}{n} + O\left(\frac{\ln n}{n^2}\right)$$

and, consequently,

$$A_1 = \frac{e^{-jx}}{n^j}\left\{1 + \frac{j(x + \ln n)}{n} + O\left(\left(\frac{\ln n}{n}\right)^2\right)\right\}. \tag{5.50}$$

Since $t_2/n \sim -1/n$ as $n \to \infty$, $A_2$ is negligible in $P\{Y_{nj} > u_n\} = A_1 - A_2$ and hence the equality (5.49) follows from (5.50). □

**Lemma 5.8.** *Let $p_{nj} = P\{Y_{nj} > n^2(x + \ln n)\}$. For any positive integer $j \geq 2$, the following relation holds as $n \to \infty$:*

$$p_{nj} - p_{n,j-1}p_{n1} = \frac{e^{-jx}}{n^{j+1}} \cdot o(1). \tag{5.51}$$

**Proof.** Relation (5.51) follows from (5.49). □

**Lemma 5.9.** *Let $x$ be a real number and $k$ and $l$ positive integers, such that $k + l \leq n$. There then exists a constant $C(x)$ such that for $u_n = n^2(x + \ln n)$, the inequality*

$$\left| P\left(\bigcap_{i=1}^{k+l}\{X_{ni} \leq u_n\}\right) - P\left(\bigcap_{i=1}^{k}\{X_{ni} \leq u_n\}\right) P\left(\bigcap_{i=k+1}^{k+l}\{X_{ni} \leq u_n\}\right) \right|$$

$$\leq C(x) \min\{k, l\} \frac{1}{n^2} \leq \frac{C(x)}{n}$$

*holds, that is, the condition $D(u_n)$ is satisfied, where this condition is defined in Chapter 3, Section 3.2 of Leadbetter et al. (1983).*

**Proof.** Let

$$\Delta_n(k, 1) = P\left(\bigcap_{j=1}^{k+1}\{X_{nj} \leq u_n\}\right) - P\left(\bigcap_{j=1}^{k}\{X_{nj} \leq u_n\}\right) \cdot P\{X_{n,k+1} \leq u_n\},$$

$$D_j = \{X_{nj} > u_n\}, \quad j = 1, 2, \ldots, k \quad \text{and} \quad A = \{X_{n,k+1} > u_n\}.$$

We then have

$$\Delta_n(k, 1) = P(D_1^c \ldots D_k^c A^c) - P(D_1^c \ldots D_k^c) P(A^c)$$
$$= 1 - P(D_1 \cup \cdots \cup D_k \cup A) - (1 - P(D_1 \cup \cdots \cup D_k))(1 - P(A))$$
$$= P((D_1 \cup \cdots \cup D_k) \cap A) - P(D_1 \cup \cdots \cup D_k) P(A)$$
$$= P(D_1 A \cup \cdots \cup D_k A) - P(D_1 \cup \cdots \cup D_k) P(A).$$

Using the inclusion–exclusion principle, we obtain that

$$\Delta_n(k, 1) = \sum_{m=2}^{k+1} (-1)^m \binom{k}{m-1} (p_{nm} - p_{n,m-1} p_{n1}). \tag{5.52}$$

Using (5.51) and (5.52), we get the statement of the lemma, first for integers $k$ and $l = 1$, then for arbitrary $k$ and $l$, where $k + l \leq n$. □

**Lemma 5.10.** *For $u_n = u_n(x) = n^2(x + \ln n)$, the condition $D'(u_n)$ is satisfied (this condition is defined in Chapter 3, Section 3.4 of Leadbetter et al. (1983)):*

$$\lim_{k \to \infty} \limsup_{n \to \infty} n \cdot \sum_{j=2}^{[n/k]} P\{X_{n1} > u_n, X_{nj} > u_n\} = 0. \tag{5.53}$$

**Proof.** It follows from (2.9) that $P\{X_{n1} > u_n, X_{nj} > u_n\} = P\{Y_{n2} > u_n\}$ holds for every $j \geq 2$. Hence, as $n \to \infty$, we get

$$n \sum_{j=2}^{[n/k]} P\{X_{n1} > u_n, X_{nj} > u_n\} = n\left(\left[\frac{n}{k}\right] - 1\right)\frac{\mathrm{e}^{-2x}}{n^2}(1 + o(1))$$

$$= \frac{\mathrm{e}^{-2x}}{k}(1 + o(1))$$

and, consequently, (5.53) holds.

Now, Theorem 3.3 follows from Lemma 5.6, Lemma 5.9, Lemma 5.10 and Theorem 5.3.1 from Leadbetter, Lindgren and Rootzén (1983). □

## Acknowledgements

The work is supported by the Ministry of Science and Environmental Protection of the Republic of Serbia, Grant No. 144032 and Grant No. 149041.

The author would like to thank the anonymous referee and Editors for useful comments and suggestions.